# EFFICIENT ESTIMATION OF A SEMIPARAMETRIC PARTIALLY LINEAR VARYING COEFFICIENT MODEL

BY IBRAHIM AHMAD, SITTISAK LEELAHANON AND QI LI[1]

*University of Central Florida, Thammasat University
and Texas A&M University*

In this paper we propose a general series method to estimate a semiparametric partially linear varying coefficient model. We establish the consistency and $\sqrt{n}$-normality property of the estimator of the finite-dimensional parameters of the model. We further show that, when the error is conditionally homoskedastic, this estimator is semiparametrically efficient in the sense that the inverse of the asymptotic variance of the estimator of the finite-dimensional parameter reaches the semiparametric efficiency bound of this model. A small-scale simulation is reported to examine the finite sample performance of the proposed estimator, and an empirical application is presented to illustrate the usefulness of the proposed method in practice. We also discuss how to obtain an efficient estimation result when the error is conditional heteroskedastic.

**1. Introduction.** Semiparametric and nonparametric estimation techniques have attracted much attention among statisticians and econometricians. One popular semiparametric specification is a partially linear model as considered by Robinson (1988), Speckman (1988) and Stock (1989), among others, via

(1) $$Y_i = v_i'\gamma + \delta(z_i) + u_i, \qquad i = 1, \ldots, n,$$

where the prime denotes transpose, $v_i'\gamma$ is the parametric component and $\delta(z_i)$ is an unknown function and, therefore, is the nonparametric component of the model; see Green and Silverman (1994), Härdle, Liang and Gao (2000) and the references therein for more detailed discussion of this model. This

Received March 2002; revised February 2004.
[1]Supported in part by the Private Enterprise Research Center, Texas A&M University.
*AMS 200 subject classification.* 62G08.

model can be generalized to the following semiparametric varying coefficient model:

$$Y_i = v_i'\gamma(z_i) + \delta(z_i) + u_i, \qquad i = 1, \ldots, n, \tag{2}$$

where $\gamma(z)$ is a vector of unknown smooth functions of $z$. Define $x_i = (1, v_i')'$ and $\beta(z) = (\delta(z), \gamma(z)')'$. Then (2) can be written more compactly as

$$Y_i = x_i'\beta(z_i) + u_i, \qquad i = 1, \ldots, n. \tag{3}$$

The varying coefficient model is an appropriate setting, for example, in the framework of a cross-sectional production function where $v_i = (\text{Labor}_i, \text{Capital}_i)'$ represents the firm's labor and capital inputs, and $z_i = \text{R\&D}_i$ is the firm's research and development expenditure. The varying coefficient model suggests that the labor and capital input coefficients may vary directly with the firm's R&D input, so the marginal productivity of labor and capital depend on the firm's R&D values. While the partially linear model (1) only allows the R&D variable to have a neutral effect on the production function, that is, it only shifts the level of the production frontier, it does not affect the labor and/or capital marginal productivity. Li, Huang, Li and Fu (2002) use the nonparametric kernel method to estimate the semiparametric varying coefficient model (2) and apply the method to China's nonmetal mineral manufacturing industry data; their results show that the semiparametric varying coefficient model (2) is more appropriate than either a parametric linear model or a semiparametric partially linear model for studying the production efficiency in China's nonmetal mineral manufacturing industry.

The time-series smooth transition autoregressive (STAR) model is another example of the varying coefficient model. It is given by $y_t = x_t'\beta(y_{t-d}) + u_t$, where $\beta(y_{t-d})$ is a vector of bounded functions; see Chen and Tsay (1993) and Hastie and Tibshirani (1993). They consider an autoregressive model of the form $y_t = f_1(y_{t-d})y_{t-1} + f_2(y_{t-d})y_{t-2} + \cdots + f_p(y_{t-d})y_{t-p} + u_t$, where the functional forms of the $f_j(\cdot)$'s ($j = 1, \ldots, p$) are not specified. Chen and Tsay (1993) and Hastie and Tibshirani (1993) discuss the identification of $f_j(\cdot)$ and suggest some recursive algorithms to estimate the unknown function $f_j(\cdot)$. More recent work on varying coefficient models can be found in Carroll, Fan, Gijbels and Wand (1997) and Fan and Zhang (1999), who propose a two-step procedure to accommodate varying degrees of smoothness among coefficient functions. See also Hoover, Rice, Wu and Yang (1998), Xia and Li (1999), Cai, Fan and Yao (2000), Cai, Fan and Li (2000), Fan and Huang (2002) and Zhang, Lee and Song (2002) on efficient estimation and inference of semiparametric varying coefficient models by using the local polynomial method and Fan, Yao and Cai (2003) on adaptive estimation of varying coefficient models.

The semiparametric varying coefficient model has the advantage that it allows more flexibility in functional forms than a parametric linear model or a semiparametric partially linear model, and, at the same time, it avoids much of the "curse of dimensionality" problem, as the nonparametric functions are restricted only to part of the variable $z$. However, when some of the $\beta$ coefficients are indeed constants, one should model them as constants and, in this way, one can obtain more efficient estimation results by incorporating this information. Consider again the production function example: if one further separates the capital into liquid capital and fixed capital, it is likely that the level of R&D will affect the marginal productivity of fixed capital, but not that of liquid capital. This gives rise to a partially linear varying coefficient model as follows:

$$(4) \qquad Y_i = w_i'\gamma + x_i'\beta(z_i) + u_i, \qquad i = 1, \ldots, n,$$

where $w_i$ is a vector of variables whose coefficient $\gamma$ is a vector of constant parameters, and say, $w$ is the firm's liquid capital in the above production example.

In this paper we propose to estimate the partially linear varying coefficient model (4) using the general series method, such as spline or power series. We show that the series method leads to efficient estimation for the finite-dimensional parameter $\gamma$ under the conditional heteroskedastic error condition. Recently, Fan and Huang (2002) suggested using the kernel-based profile likelihood approach to estimate a partially varying coefficient model [this paper was brought to our attention after the first submission of our paper], and they show that their approach also leads to efficient estimation of the finite-dimensional parameter $\gamma$ when the error is conditional homoskedastic. In this paper we also argue that the efficient estimation result of the series-based method can be extended to the conditional heteroskedastic error case in a straightforward way. It is more difficult to obtain efficient estimation results using the kernel-based method when the error is conditional heteroskedastic. Moreover, the series estimators have well-defined meanings as estimating the best approximation function for the unknown conditional mean regression function even when the model is misspecified. The payoff of using the general series estimation methods is that it is difficult to establish the asymptotic normality result for the nonparametric components under optimal smoothings (i.e., balance the squared bias and variance terms). Thus, the series method should be viewed as a complement to the kernel method in estimating a partially linear varying coefficient model.

**2. Estimation.** Consider the following partially linear varying coefficient model:

$$(5) \qquad Y_i = w_i'\gamma + x_i'\beta(z_i) + u_i, \qquad i = 1, \ldots, n,$$

where $w_i$ is a $q \times 1$ vector of random variables, $\gamma$ is a $q \times 1$ vector of unknown parameters, $x_i$ is of dimension $d \times 1$, $z_i = (z_{i1}, \ldots, z_{ir})$ is of dimension $r$, $\beta(\cdot) = (\beta_1(\cdot), \ldots, \beta_d(\cdot))'$ is a $d \times 1$ vector of unknown varying coefficient functions, and $u_i$ is an error term satisfying $E(u_i|w_i, x_i, z_i) = 0$.

With the series estimation method, for $l = 1, \ldots, d$, we approximate the varying coefficient function $\beta_l(z)$ by $p_l^{k_l}(z)'\alpha_l^{k_l}$, a linear combination of $k_l$ base functions, where $p_l^{k_l}(z) = [p_{l1}(z), \ldots, p_{lk_l}(z)]'$ is a $k_l \times 1$ vector of base functions and $\alpha_l^{k_l} = (\alpha_{l1}, \ldots, \alpha_{lk_l})'$ is a $k_l \times 1$ vector of unknown parameters. The approximation functions $p_l^{k_l}(z)$ have the property that, as $k_l$ grows, there is a linear combination of $p_l^{k_l}(z)$ that can approximate any smooth function $\beta_l(z)$ arbitrarily well in the sense that the approximation mean square error can be made arbitrarily small.

Define the $K \times 1$ matrices $p^K(x_i, z_i) = (x_{i1}p_1^{k_1}(z_i)', \ldots, x_{id}p_d^{k_d}(z_i)')'$ and $\alpha = (\alpha_1^{k_1\prime}, \ldots, \alpha_d^{k_d\prime})'$, where $K = \sum_{l=1}^d k_l$. Thus, we use a linear combination of $K$ functions, $p^K(x_i, z_i)'\alpha$, to approximate $x_i'\beta(z_i)$. Hence, we can rewrite (5) as

$$
\begin{aligned}
Y_i &= w_i'\gamma + p^K(x_i, z_i)'\alpha + (x_i'\beta(z_i) - p^K(x_i, z_i)'\alpha) + u_i \\
&= w_i'\gamma + p_i^K(x_i, z_i)'\alpha + \text{error}_i,
\end{aligned}
\tag{6}
$$

where the definition of $\text{error}_i$ should be apparent.

We introduce some matrix notation. Let $Y = (Y_1, \ldots, Y_n)'$, $u = (u_1, \ldots, u_n)'$, $W = (w_1, \ldots, w_n)'$, $G = (x_1'\beta(z_1), \ldots, x_n'\beta(z_n))'$ and $P = (p^K(x_1, z_1), \ldots, p^K(x_n, z_n))'$. Hence, model (6) can be written in matrix notation as

$$Y = W\gamma + P\alpha + \text{error}. \tag{7}$$

Let $\hat{\gamma}$ and $\hat{\alpha}$ denote the least squares estimators of $\gamma$ and $\alpha$ obtained by regressing $Y$ on $(W, P)$ from (7). Then we estimate $\beta_l(z)$ by $\hat{\beta}_l(z) \stackrel{\text{def}}{=} p_l^{k_l}(z)'\hat{\alpha}_l$ ($l = 1, \ldots, d$). We will establish the $\sqrt{n}$-normality result for $\hat{\gamma}$ and derive the rate of convergence for $\hat{\beta}_l(z)$.

We present an alternative form for $\hat{\gamma}$ and $\hat{\alpha}$ that is convenient for the asymptotic analysis given below. In matrix form, (5) can be written as

$$Y = W\gamma + G + u. \tag{8}$$

Define $M = P(P'P)^- P'$, where $(\cdot)^-$ denotes any symmetric generalized inverse of $(\cdot)$. [Under the assumptions given in this paper, $P'P$ is nonsingular with probability one. In finite sample applications, if $P'P$ is singular, one can remove the redundant regressors to make $P'P$ nonsingular.] For an $n \times m$ matrix $A$, define $\widetilde{A} = MA$. Then premultiplying (8) by $M$ leads to

$$\widetilde{Y} = \widetilde{W}\gamma + \widetilde{G} + \tilde{u}. \tag{9}$$

Subtracting (9) from (8) yields

$$Y - \widetilde{Y} = (W - \widetilde{W})\gamma + (G - \widetilde{G}) + u - \tilde{u}. \tag{10}$$

$\hat{\gamma}$ can also be obtained as the least squares regression of $Y - \widetilde{Y}$ on $W - \widetilde{W}$, that is,

$$\hat{\gamma} = [(W - \widetilde{W})'(W - \widetilde{W})]^{-}(W - \widetilde{W})'(Y - \widetilde{Y}). \tag{11}$$

And $\hat{\alpha}$ can be obtained from (7) with $\gamma$ being replaced by $\hat{\gamma}$,

$$\hat{\alpha} = (P'P)^{-} P'(Y - W\hat{\gamma}), \tag{12}$$

from which we obtain $\hat{\beta}_l(z) = p_l^{k_l}(z)' \hat{\alpha}_l^{k_l}$, $l = 1, \ldots, d$.

Under the assumptions given below, both $(W - \widetilde{W})'(W - \widetilde{W})$ and $P'P$ are asymptotically nonsingular. Hence, $\hat{\gamma}$ and $\hat{\alpha}$ given in (11) and (12) are well defined and they are numerically identical to the least squares estimator obtained by regressing $Y$ on $(W, P)$.

Next we give a definition and some assumptions that are used to derive the main results of this paper.

DEFINITION 2.1. $g(x, z)$ is said to belong to the varying coefficient class of functions $\mathcal{G}$ if:

(i) $g(x, z) = x'h(z) \equiv \sum_{l=1}^{d} x_l h_l(z)$ for some continuous functions $h_l(z)$, where $h(z) = (h_1(z), \ldots, h_d(z))'$.
(ii) $\sum_{l=1}^{d} E[x_{il}^2 h_l(z_i)^2] < \infty$, where $x_l$ ($x_{il}$) is the $l$th component of $x$ ($x_i$).

For any function $f(x, z)$, let $E_{\mathcal{G}}[f(x, z)]$ denote the projection of $f(x, z)$ onto the varying coefficient functional space $\mathcal{G}$ (under the $L_2$-norm). That is, $E_{\mathcal{G}}[f(x, z)]$ is an element that belongs to $\mathcal{G}$ and it is the closest function to $f(x, z)$ among all the functions in $\mathcal{G}$. More specifically ($x_l$ is the $l$th component of $x$, $l = 1, \ldots, d$),

$$E\{(f(x,z) - E_{\mathcal{G}}[f(x,z)])(f(x,z) - E_{\mathcal{G}}[f(x,z)])'\}$$
$$= \inf_{\sum_l x_l h_l(z) \in \mathcal{G}} E\left\{\left(f(x,z) - \sum_{l=1}^{d} x_l h_l(z)\right)\left(f(x,z) - \sum_{l=1}^{d} x_l h_l(z)\right)'\right\}. \tag{13}$$

Thus,

$$E[(f(x,z) - E_{\mathcal{G}}[f(x,z)])(f(x,z) - E_{\mathcal{G}}[f(x,z)])']$$
$$\leq E\left[\left(f(x,z) - \sum_{l=1}^{d} x_l h_l(z)\right)\left(f(x,z) - \sum_{l=1}^{d} x_l h_l(z)\right)'\right], \tag{14}$$

for all $g(x, z) = \sum_{l=1}^{d} x_l h_l(z) \in \mathcal{G}$. Here for square matrices $A$ and $B$, $A \leq B$ means that $A - B$ is negative semidefinite.

Define $\theta(x, z) = E[w|x, z]$ and $m(x, z) = E_{\mathcal{G}}[\theta(x, z)]$. The following assumptions will be used to establish the asymptotic distribution of $\hat{\gamma}$ and the convergence rates of $\hat{\beta}(z)$.

ASSUMPTION 2.1. (i) $(Y_i, w_i, x_i, z_i)_{i=1}^n$ are independent and identically distributed as $(Y_1, w_1, x_1, z_1)$ and the support of $(w_1, x_1, z_1)$ is a compact subset of $\mathcal{R}^{q+d+r}$; (ii) both $\theta(x_1, z_1)$ and $\text{var}[Y_1|w_1, x_1, z_1]$ are bounded functions on the support of $(w_1, x_1, z_1)$.

ASSUMPTION 2.2. (i) For every $K$ there is a nonsingular matrix $B$ such that for $P^K(x,z) = B p^K(x,z)$ the smallest eigenvalue of $E[P^K(x_i, z_i) P^K(x_i, z_i)']$ is bounded away from zero uniformly in $K$; (ii) there is a sequence of constants $\zeta_0(K)$ satisfying $\sup_{(x,z) \in \mathcal{S}} \|P^K(x,z)\| \leq \zeta_0(K)$ and $K = K_n$ such that $(\zeta_0(K))^2 K/n \to 0$ as $n \to \infty$, where $\mathcal{S}$ is the support of $(x_1, z_1)$, and for a matrix $A$, $\|A\| = [\text{tr}(A'A)]^{1/2}$ denotes the Euclidean norm of $A$.

ASSUMPTION 2.3. (i) For $f(x,z) = \sum_{l=1}^d x_l \beta_l(z)$ or $f(x,z) = m_j(x,z)$ $(j=1,\ldots,q)$, there exist some $\delta_l > 0$ $(l=1,\ldots,d)$, $\alpha_f = \alpha_{fK} = (\alpha_1^{k_1\prime}, \ldots, \alpha_d^{k_d\prime})'$, such that $\sup_{(x,z) \in \mathcal{S}} |f(x,z) - P^K(x,z)' \alpha_f| = O(\sum_{l=1}^d k_l^{-\delta_l})$; (ii) for $\min\{k_1, \ldots, k_d\} \to \infty$, $\sqrt{n}(\sum_{l=1}^d k_l^{-2\delta_l}) \to 0$ as $n \to \infty$.

Assumption 2.1 is a standard assumption being used on series estimation methods. Assumption 2.2 usually implies that the density function of $(x, z)$ needs to be bounded below by a positive constant. Assumption 2.3 says that there exist some $\delta_l > 0$ $(l = 1, \ldots, d)$ such that the uniform approximation error to the function shrinks at the rate $\sum_{l=1}^d k_l^{-\delta_l}$. Assumptions 2.2 and 2.3 are not the easiest conditions, but it is known that many series functions satisfy these conditions, for example, power series and splines.

Under the above assumptions, we can state our main theorem.

THEOREM 2.1. *Define $\varepsilon_i = w_i - m(x_i, z_i)$, where $m(x_i, z_i) = E_\mathcal{G}(w_i)$, and assume that $\Phi \equiv E[\varepsilon_i \varepsilon_i']$ is positive definite. Then under Assumptions 2.1–2.3 we have:*

(i) $\sqrt{n}(\hat{\gamma} - \gamma) \to N(0, \Sigma)$ *in distribution, where* $\Sigma = \Phi^{-1} \Omega \Phi^{-1}$, $\Omega = E[\sigma^2(w_i, x_i, z_i) \varepsilon_i \varepsilon_i']$ *and* $\sigma^2(w_i, x_i, z_i) = E[u_i^2 | w_i, x_i, z_i]$.

(ii) *A consistent estimator of $\Sigma$ is given by $\hat{\Sigma} = \hat{\Phi}^{-1} \hat{\Omega} \hat{\Phi}^{-1}$, where $\hat{\Phi} = n^{-1} \sum_{i=1}^n (w_i - \tilde{w}_i)(w_i - \tilde{w}_i)'$, $\hat{\Omega} = n^{-1} \sum_{i=1}^n \hat{u}_i^2 (w_i - \tilde{w}_i)(w_i - \tilde{w}_i)'$, $\tilde{w}_i$ is the ith row of $\widetilde{W}$ and $\hat{u}_i = Y_i - w_i' \hat{\gamma} - p^K(x_i, z_i)' \hat{\alpha}$.*

The proof of Theorem 2.1 is given in the Appendix. [One may prove Theorem 2.1 based on the general result of Shen (1997) and Ai and Chen (2003) which requires one to establish stochastic equicontinuity of the objective function. However, for the specific partially linear varying semiparametric model, it is easier to use a direct proof as given in the Appendix.]

Under the conditional homoskedastic error assumption $E[u_i^2|w_i, x_i, z_i] = E(u_i^2) = \sigma^2$, the estimator $\hat{\gamma}$ is semiparametric efficient in the sense that the inverse of the asymptotic variance of $\sqrt{n}(\hat{\gamma} - \gamma)$ equals the semiparametric efficiency bound. From the result of Chamberlain (1992) [the concept of semiparametric efficient bound we use here is discussed in Chamberlain (1992), which gives the lower bound for the asymptotic variance of an (regular) estimator satisfying some conditional moment conditions; see also Bickel, Klaassen, Ritov and Wellner (1993) for a more general treatment of efficient and adaptive inference in semiparametric models], the semiparametric efficiency bound for the inverse of the asymptotic variance of an estimator of $\gamma$ is

$$(15) \quad J_0 = \inf_{g \in \mathcal{G}} E[(w_i - g(x_i, z_i))(\text{var}[u_i|w_i, x_i, z_i])^{-1}(w_i - g(x_i, z_i))'].$$

Under the conditional homoskedastic error assumption $\text{var}[u_i|w_i, x_i, z_i] = \sigma^2$, then (15) can be rewritten as $(m(x_i, z_i) = E_{\mathcal{G}}(w_i))$

$$
\begin{aligned}
J_0 &= \frac{1}{\sigma^2} \inf_{g \in \mathcal{G}} E[(w_i - g(x_i, z_i))(w_i - g(x_i, z_i))'] \\
&= \frac{1}{\sigma^2} E[(w_i - m(x_i, z_i))(w_i - m(x_i, z_i))'] \\
&= \frac{1}{\sigma^2} E[\varepsilon_i \varepsilon_i'] = \frac{\Phi}{\sigma^2}.
\end{aligned}
\tag{16}
$$

Note that the inverse of (16) coincides with $\Sigma = \sigma^2 \Phi^{-1}$, the asymptotic variance of $\sqrt{n}(\hat{\gamma} - \gamma)$ when the error is conditional homoskedastic. Hence, $\Sigma^{-1} = J_0$ and $\hat{\gamma}$ is a semiparametrically efficient estimator under the conditional homoskedastic error assumption.

The next theorem gives the convergence rate of $\hat{\beta}_l(z) = p_l^{k_l}(z)\hat{\alpha}_l^{k_l}$ to $\beta_l(z)$ for $l = 1, \ldots, d$.

THEOREM 2.2. *Under Assumptions 2.1–2.3, let $\mathcal{S}_z$ denote the support of $z_i$. Then we have, for $l = 1, \ldots, d$:*

(i) $\sup_{z \in \mathcal{S}_z} |\hat{\beta}_l(z) - \beta_l(z)| = O_p(\zeta_0(K)(\sqrt{K}/\sqrt{n} + \sum_{l=1}^d k_l^{-\delta_l}))$.

(ii) $\frac{1}{n} \sum_{i=1}^n (\hat{\beta}_l(z) - \beta_l(z))^2 = O_p(K/n + \sum_{l=1}^d k_l^{-2\delta_l})$.

(iii) $\int (\hat{\beta}_l(z) - \beta_l(z))^2 \, dF_z(z) = O_p(K/n + \sum_{l=1}^d k_l^{-2\delta_l})$, *where $F_z$ is the cumulative distribution function of $z_i$.*

The proof of Theorem 2.2 is given in the Appendix.

Newey (1997) gives some primitive conditions for power series and $B$-splines such that the Assumptions 2.1–2.3 hold. We state them here for the readers' convenience.

ASSUMPTION 2.4. (i) The support of $(x_i, z_i)$ is a Cartesian product of compact connected intervals on which $(x_i, z_i)$ has an absolutely continuous probability density function that is bounded above by a positive constant and bounded away from zero; (ii) for $l = 1, \ldots, d$, $f_l(x, z)$ is continuously differentiable of order $c_l$ on the support $\mathcal{S}$, where $f_l(x, z) = x_l \beta_l(z)$ or $f_l(x, z) = m_l(x, z)$.

ASSUMPTION 2.5. The support of $(x_i, z_i)$ is $[-1, 1]^{d+r}$.

Suppose that a smooth function $\eta(z)$ ($z \in R^r$) is continuously differentiable of order $c$. It is well established that the approximation error by using power series or $B$-splines is of the order of $O(K^{-c/r})$; see Lorentz (1966), Andrews (1991), Newey (1997) and Huang (1998). Therefore, Assumption 2.3(i) holds for power series and $B$-splines if $\sqrt{n}(\sum_{l=1}^{d} k_l^{-c_l/r}) = o(1)$ (i.e., $\delta_l = c_l/r$). Newey (1997) shows that, for power series or splines, Assumption 2.4 implies that the smallest eigenvalue of $E[P^K(x_i)P^K(x_i)']$ is bounded for all $K$. Also, Assumptions 2.4 and 2.5 imply that Assumptions 2.2 and 2.3 hold for $B$-splines with $\zeta_0(K) = O(\sqrt{K})$. Hence, we have the following results for regression splines.

THEOREM 2.3. *For splines, if Assumptions* 2.1, 2.4 *and* 2.5 *are satisfied, and $k_l^2/n \to 0$ as $n \to \infty$ for $l = 1, \ldots, d$, then:*

(i) *The conclusion of Theorem* 2.1 *holds.*
(ii) *The conclusion of Theorem* 2.2 *holds with $\sqrt{K}$ replacing $\zeta_0(K)$.*

Theorem 2.2 only gives the rate of convergence of the series estimator for the varying coefficient function $\beta(z)$. As we mentioned in the Introduction, it is difficult to obtain asymptotic normality results for the series estimator of $\beta(z)$ under optimal smoothings. The reason is that the asymptotic bias of the series estimator is unknown in general. Recently, Zhou, Shen and Wolfe (1998) have obtained an asymptotic bias for univariate spline regression functions that belong to $C^p$ (i.e., the regression functions have continuous $p$th derivatives) under somewhat stringent conditions such as the knots are asymptotically equally-spaced, and the degree of the spline $m$ is equal to $p - 1$. See Huang (2003) for a more detailed discussion on the difficulty of obtaining the asymptotic bias for general cases with splines. Alternatively, one may choose to undersmooth the data. In this case the bias is asymptotically negligible. Huang (2003) has obtained the asymptotic distribution of spline estimators under quite general conditions (provided the data are slightly undersmoothed). Huang, Wu and Zhou (2002, 2004) have further provided asymptotic distribution results for spline estimation of a varying coefficient model. Their results can be directly applied to obtain the asymptotic distribution of $\hat{\beta}(z)$ in a partially linear varying coefficient model. This is because $\hat{\gamma} - \gamma = O_p(n^{-1/2})$, which converges to zero faster than any nonparametric estimation convergence rate. Therefore, $\hat{\beta}(z)$ has the same asymptotic distribution whether one uses the estimator $\hat{\gamma}$ or the true $\gamma$, the latter becomes a varying coefficient model (when $\gamma$ is unknown) and the results of Huang, Wu and Zhou (2002, 2004) apply.

**3. Monte Carlo simulations.** In this section we report some simulation results to examine the finite sample performance of our proposed estimator, and also compare it with the kernel-based profile likelihood estimator suggested by Fan and Huang (2002). We first consider the following data generating process (DGP):

(17) $\quad \text{DGP1} : y_i = 1 + 0.5 w_i + x_i \beta_1(z_i) + u_i, \qquad i = 1, \ldots, n,$

where

(18) $\quad \beta_1(z_i) = 1 + (24 z_i)^3 \exp(-24 z_i)$

is taken from Hart (1997), $\beta_0 = 1$ and $\gamma = 0.5$. The error $u_i$'s are i.i.d. normal with mean 0 and variance 0.25, $z_i$ is generated by the i.i.d. uniform$[0, 2]$ distribution, $w_i = v_{1i} + 2 v_{3i}$ and $x_i = v_{2i} + v_{3i}$, where $v_{ji}$, $j = 1, 2, 3$, are i.i.d. uniform$[0, 2]$.

We also consider a second data generating process:

(19) $\quad \text{DGP2} : y_i = 4 + 0.5 w_i + x_{i1} \beta_1(z_i) + x_{i2} \beta_2(z_i) + u_i, \qquad i = 1, \ldots, n,$

where $\beta_1(z_i)$ is the same as in DGP1, $\beta_2(z_i) = z_i + \sin(z_i)$, $z_i$ is i.i.d. uniform$[0, 2]$, $u_i$ is i.i.d. normal with mean 0 and variance 0.25, $w_i = v_{1i} + 2 v_{3i}$, $x_{1i} = v_{2i} + v_{3i}$, and $x_{2i} = v_{4i} + 0.5 v_{3i}$, where $v_{ji}$ ($j = 1, 2, 3, 4$) are i.i.d. uniform$[0, 2]$.

The sample sizes are $n = 100$ and $n = 200$, and the number of replications is 5000 for all cases. We compare the estimated mean squared error (MSE) of $\hat{\gamma}$ defined by $MSE(\hat{\gamma}) = \frac{1}{5000} \sum_{j=1}^{5000} (\hat{\gamma}_j - \gamma)^2$, and estimated mean average squared error (MASE) of $\hat{\beta}_l(\cdot)$ defined by $MASE(\hat{\beta}_l(\cdot)) = \frac{1}{5000} \sum_{j=1}^{5000} [\frac{1}{n} \sum_{i=1}^{n} (\hat{\beta}_{l,j}(z_i) - \beta_l(z_i))^2$ ($l = 1$ for DGP1, $l = 1, 2$ for DGP2), where $\hat{\gamma}_j$ and $\hat{\beta}_{l,j}(z_i)$ are, respectively, the estimates of $\gamma$ and $\beta_l(z_i)$ from the $j$th replication based on one of the two methods: the $B$-spline method and the kernel-based profile likelihood method. We use a univariate cubic $B$-spline basis function defined by

(20) $\quad B(z | t_0, \ldots, t_4) = \frac{1}{3!} \sum_{j=0}^{4} (-1)^j \binom{4}{j} [\max(0, z - t_j)]^3,$

where $t_0, \ldots, t_4$ are the evenly-spaced design knots. The kernel estimator of $\gamma$ is discussed at the end of Section 2. The number of terms $K$ in series estimation and the smoothing parameter $h$ in kernel estimation are both selected by leave-one-out least squares cross-validation. As discussed in Bickel and Kwon (2002), the estimation of the parametric component does not very sensitively depend on the choice of smoothing parameters, as long as the selected smoothing parameters do not create excessive bias in the estimation of the nonparametric components. In this regard, the cross-validation method usually performs well. (Other data driven methods in selecting $K$ in series estimation include the following: the generalized cross-validation criterion [Craven and Wahba (1979) and Li (1987)] and Mallows' $C_p$ criterion [Mallows (1973)].)

TABLE 1
$MSE(\hat{\gamma})$ *by spline and kernel methods*

|  |  | DGP1 | | DGP2 | |
|---|---|---|---|---|---|
|  |  | $n = 100$ | $n = 200$ | $n = 100$ | $n = 200$ |
| Cubic $B$-spline | $MSE(\hat{\gamma})$ | 0.00278 | 0.00133 | 0.00357 | 0.00153 |
| Profile likelihood | $MSE(\hat{\gamma})$ | 0.00315 | 0.00145 | 0.00443 | 0.00178 |

The simulation result is presented in Table 1. From Table 1, first we observe that as the sample size doubles, the estimated MSE for all three different estimators reduces to about half of the original values; this is consistent with the fact that all of them are $\sqrt{n}$-consistent estimators of $\gamma$. Second we observe that the $B$-spline method gives slightly smaller estimated MSE of $\hat{\gamma}$ for both DGPs. Under the conditional homoskedastic error condition, both methods are semiparametrically efficient. Therefore, they have the same asymptotic efficiency. The results in Table 1 may reflect small sample differences of the two estimation methods for the chosen data generating processes (DGP). It is possible that for some other DGPs the kernel method may have better small sample performance. In fact, a few simulated examples cannot differentiate the finite sample performance of the two methods.

Table 2 reports $MASE(\hat{\beta}(z))$ for the spline and the profile likelihood methods. The spline and the kernel methods give similar estimation results for $MASE(\hat{\beta}(z))$ for both DGPs.

The results of Tables 1 and 2 are based on the least-squares cross-validation selection of $K$ (for spline) and $h$ (for the profile likelihood method). To examine whether our findings only reflect a particular way of selecting the smoothing parameters (the cross-validation method), we also compute the $MSE(\hat{\gamma})$ and $MASE(\hat{\beta}(\cdot))$ for a range of different values of $K$ and $h$ without leave-one-out in the estimation. Figures 1 and 2 plot the estimation results.

In Figure 1(a) the dashed line plots the leave-one-out cross-validation function for a range of $K$ for the spline method (DGP1, $n = 100$, average over the 5,000 replications). We observe that the cross validation function

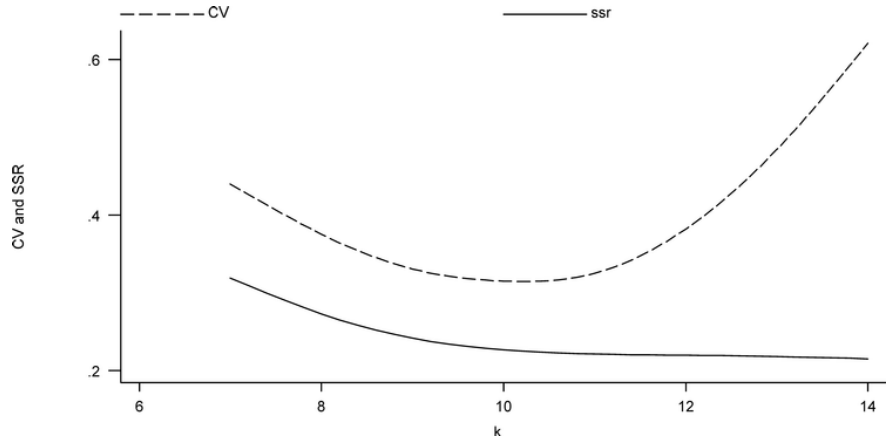

(a)

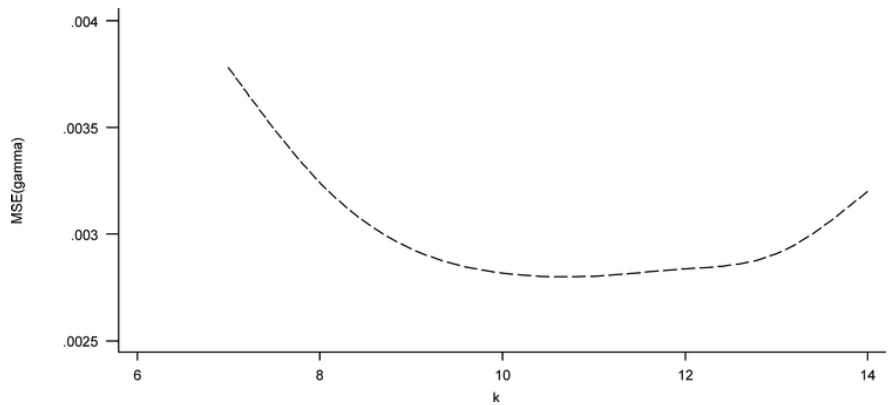

(b)

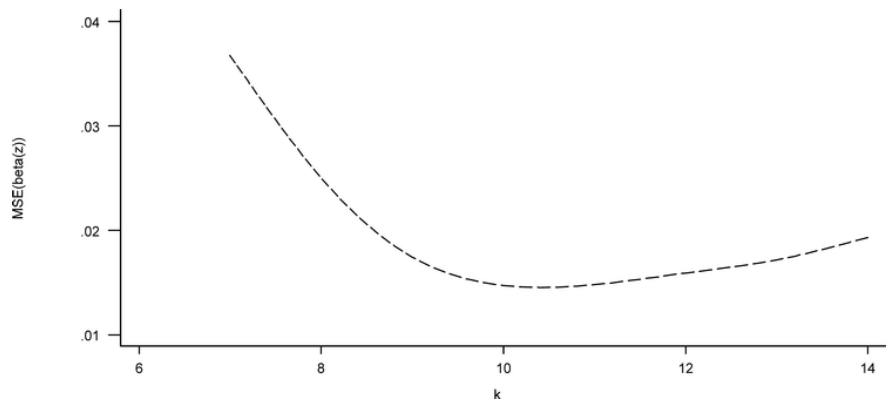

(c)

FIG. 1. (a) *CV function and SSR (spline, DGP*1*, $n = 100$).* (b) *MSE($\hat{\gamma}$) (spline, DGP*1*, $n = 100$).* (c) *MASE($\hat{\beta}(z)$) (spline, DGP*1*, $n = 100$).*

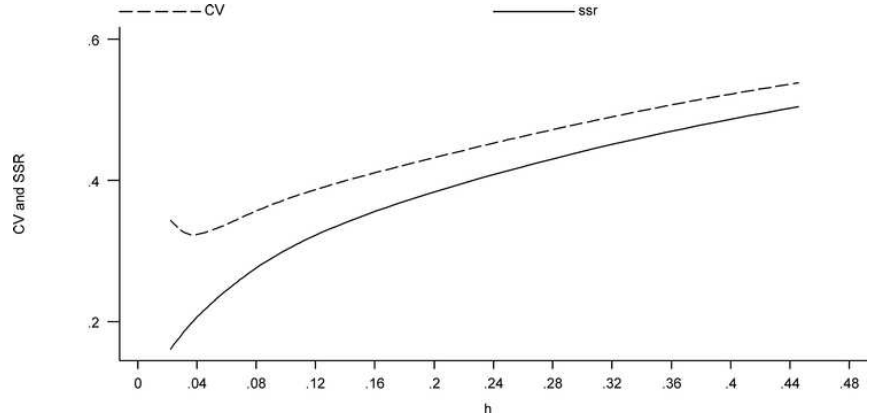

(a)

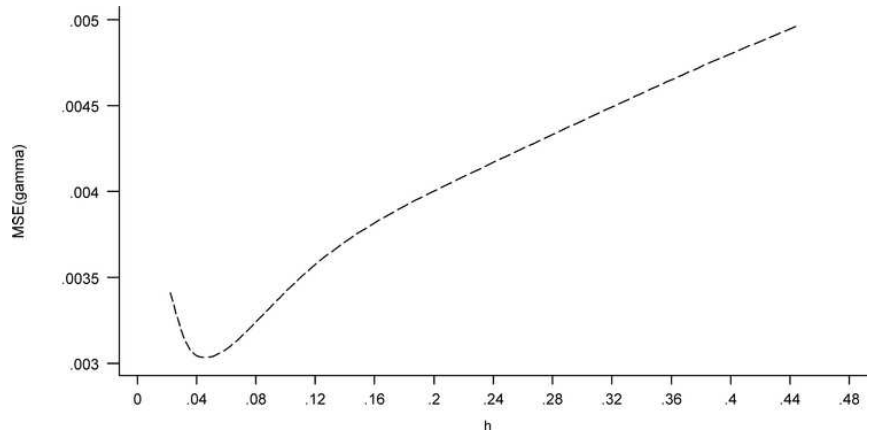

(b)

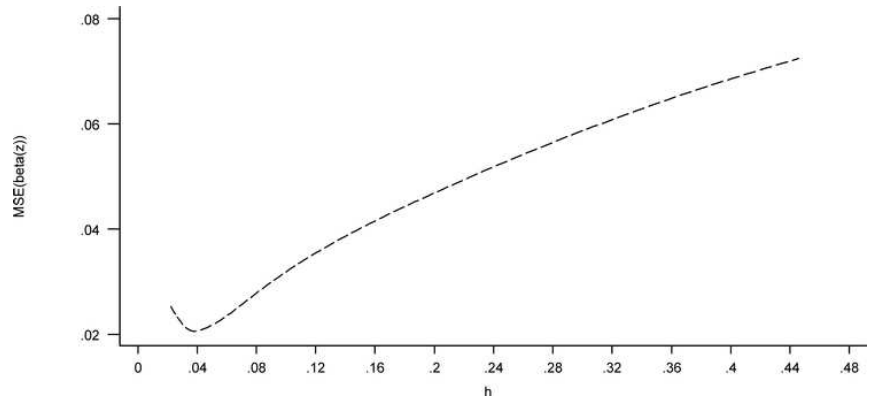

(c)

Fig. 2. (a) *CV function and SSR (kernel, DGP1, $n = 100$).* (b) *$MSE(\hat{\gamma})$ (kernel, DGP1, $n = 100$).* (c) *$MASE(\hat{\beta}(z))$ (kernel, DGP1, $n = 100$).*

is minimized around $K = 10$. The solid line in Figure 1(a) is the sum of squared residuals computed without using the leave-one-out estimator; as expected, it decreases as $K$ increases.

Figure 1(b) graphs the $MSE(\hat{\gamma})$ computed using all observations (not using the leave-one-out method). We see that $MSE(\hat{\gamma})$ takes minimum values around $K = 10$ and $K = 11$. Figure 1(c) plots the $MASE(\hat{\beta}(z))$, again computed using all observations. $MASE(\hat{\beta}(z))$ assumes minimum values around $K = 10$. The average of 5,000 cross-validation selected $K$'s is 10.42.

From Figure 1 we can see that, on average, the least squares cross-validation method performs well in selecting $K$ that is close to values of $K$ that minimize $MSE(\hat{\gamma})$ and $MASE(\hat{\beta}(z))$. Note that both Figures 1(b) and 1(c) do not use the leave-one-out estimator. Therefore, unlike the sum of squared residuals, $MSE(\hat{\gamma})$ and $MASE(\hat{\beta}(z))$ do not monotonically decrease as $K$ increases.

Figure 2 gives the corresponding cases for the profile kernel method. Figure 2(a) shows that the cross-validation function is minimized around $h = 0.04$, while the sum of squares of residuals monotonically increases with $h$.

Figures 2(b) and 2(c) show that both $MSE(\hat{\gamma})$ and $MASE(\hat{\beta}(z))$ are minimized around $h = 0.04$. Note that Figures 2(b) and 2(c) are computed using all observations (without using the leave-one-out method). Therefore, similar to the spline case, $MSE(\hat{\gamma})$ and $MASE(\hat{\beta}(z))$ do not decrease monotonically with $h$, but rather they are both minimized around the value of $h$ that minimizes the cross-validation function.

Summarizing the results of Figures 1 and 2, we find that the cross-validation method performs adequately for the simulated data. The simulation results reported in this section show that both the spline and the kernel methods can be a useful tool in estimating a partially linear varying coefficient model.

**4. An empirical application.** In this section we consider estimation of a production function in China's manufacturing industry to illustrate the application of the partially linear varying coefficient model. The data used in this paper are drawn from the Third Industrial Census of China conducted by the National Statistical Bureau of China in 1995. The Third Industrial Census of China is currently the most comprehensive industrial survey in China. To avoid heterogeneity across different industries and also to maintain enough observations in the sample for accurate semiparametric estimation, we include firms from the sector of food, soft drink and cigarette manufacturing in this study. After removing firms with missing values, the sample size we use is 877. We estimate a benchmark parametric linear model as follows:

$$\ln Y = \beta_0 + \gamma \ln w + \beta_l \ln L + \beta_k \ln K + \beta_z \ln z + u, \tag{21}$$

where $Y$ is the sales of the firm, $w$ is the liquid capital, $L$ is the labor input, $K$ is the fixed capital and $z$ is the firm's R&D (all monetary measures are in thousand RMB, the Chinese currency).

TABLE 2
$MASE(\hat{\beta}(\cdot))$ by spline and kernel methods

|  | DGP1 | | DGP2 | | | |
| --- | --- | --- | --- | --- | --- | --- |
|  | $MASE(\hat{\beta}_1(\cdot))$ | | $MASE(\hat{\beta}_1(\cdot))$ | | $MASE(\hat{\beta}_2(\cdot))$ | |
|  | $n = 100$ | $n = 200$ | $n = 100$ | $n = 200$ | $n = 100$ | $n = 200$ |
| Cubic $B$-spline | 0.0162 | 0.00764 | 0.0576 | 0.0245 | 0.0635 | 0.0326 |
| Profile likelihood | 0.0224 | 0.0110 | 0.0815 | 0.0356 | 0.0593 | 0.0318 |

The partially linear varying coefficient model is given by

$$\ln Y = \gamma \ln w + \beta_0(z) + \beta_l(z) \ln L + \beta_k(z) \ln K + u. \tag{22}$$

Here we choose liquid capital as the $w$ variable whose coefficient does not depend on the firm's R&D spending ($z$). We have given some theoretical arguments for this model specification in the Introduction; to justify this choice statistically, we test both models (21) and (22) against a more general semiparametric varying coefficient model,

$$\ln Y = \gamma(z) \ln w + \beta_0(z) + \beta_l(z) \ln L + \beta_k(z) \ln K + u. \tag{23}$$

Obviously, (23) includes (22) as special case when $\gamma(z)$ is constant for all $z$. We use quadratic and cubic splines and the number of knots is chosen by the least squares cross-validation method. The cross-validation method selected the quadratic spline. Our test for the null models (21) and (22) is based on $(RSS_0 - RSS)/RSS$, where $RSS_0$ is the residual sum of squares from the null model, and $RSS$ is from the alternative model (23). We obtain the critical values of our test based on 1,000 residual-based bootstrap procedures where we first obtain the residuals from the null model, from which we generate two point wild bootstrap errors, which in turn are used to generate bootstrap $\ln Y$'s (using the estimated null model); the bootstrap statistic is $(RSS_0^* - RSS^*)/RSS^*$, where $RSS_0^*$ is the residual sum of squares from the null model computed using the bootstrap sample and $RSS^*$ is computed from the alternative model also using the bootstrap sample. Note that the bootstrap sample is generated according to the null model. Therefore, the bootstrap statistic approximates the *null* distribution of the original test statistic even when the null hypothesis is false. When testing the parametric null model, we firmly reject the null model with a $p$-value of 0.001. For testing the partially linear varying coefficient model (22), we cannot reject this null model at conventional levels (a $p$-value of 0.162). Therefore, both economic theory and the statistical testing results support our specification (22).

The estimated value of $\gamma$ based on (22) is 0.481, with a standard error of 0.0372 (the $t$-statistic is 12.91). The goodness-of-fit $R^2$ is 0.566 $[R^2 = 1 - RSS/\sum_i (y_i - \bar{y})^2, y_i = \ln Y_i]$. The estimated varying coefficient functions are plotted in Figures 3(a) to 3(c). $\beta_0(z)$ is plotted in Figure 3(a). Figure 3(b) shows that the marginal productivity of labor $\beta_l(z)$ is a nonlinear function of $z$ (R&D). The marginal productivity of labor first increases with $z$ and then decreases as $z$ increases further. The bell shape of the curve suggests that, while modest R&D can improve labor productivity, higher R&D leads to lower labor productivity. Figure 3(c) shows that the marginal productivity of (fixed) capital is also nonlinear in $z$. It exhibits a general up trend with $z$, indicating that firms with large R&D spending yield relative higher marginal (fixed) capital productivity. These results are not surprising given that most of the firms in our sample are state-owned. It is typical in these firms that capital is scarce while labor is excessive. Thus, most of the R&D expenses are used to improve equipment performance, but not to train labor. In Figure 3(d) we graph the return to scale function $\gamma + \beta_l(z) + \beta_k(z)$. The return to scale is well below one (the constant return to scale level) for firms with small R&D, and it increases to a range between 0.8 to 0.9 for firms with large R&D expenditures. The results indicate that most of the firms in our sample exhibit decreasing returns to scale in production. It partly reflects the fact that the firms included in the survey are large firms, most of which are state-owned firms. These firms typically have a production scale larger than ideal. In particular, there are usually too many employees in these firms. It was not until several years after the survey we use in this paper, as a result of fierce competition from foreign firms and the passage of bankruptcy law in China, that the food, soft drink and cigarette sector witnessed a string of reorganizations, mergers and acquisitions. Further discussion is beyond the scope of this paper.

We have also applied the kernel profile likelihood method to this data set. The estimation results are quite similar to those obtained by the spline method. For example, the estimated $\gamma$ is 0.489 with a $t$-statistic of 13.20. The $\beta(z)$ functions all have similar shapes as those obtained by the spline method. Therefore, we do not report the kernel estimation results here.

**5. Possible extension.** In this section we briefly discuss (without providing technical details) efficient estimation of a partially varying coefficient model when the error is conditional heteroskedastic.

Theorem 2.1 holds even when the error is conditional heteroskedastic, say, $E(u_i^2|v_i) = \sigma^2(v_i)$, where $v_i = (w_i, x_i, z_i)$. However, in this case $\hat{\gamma}$ is not

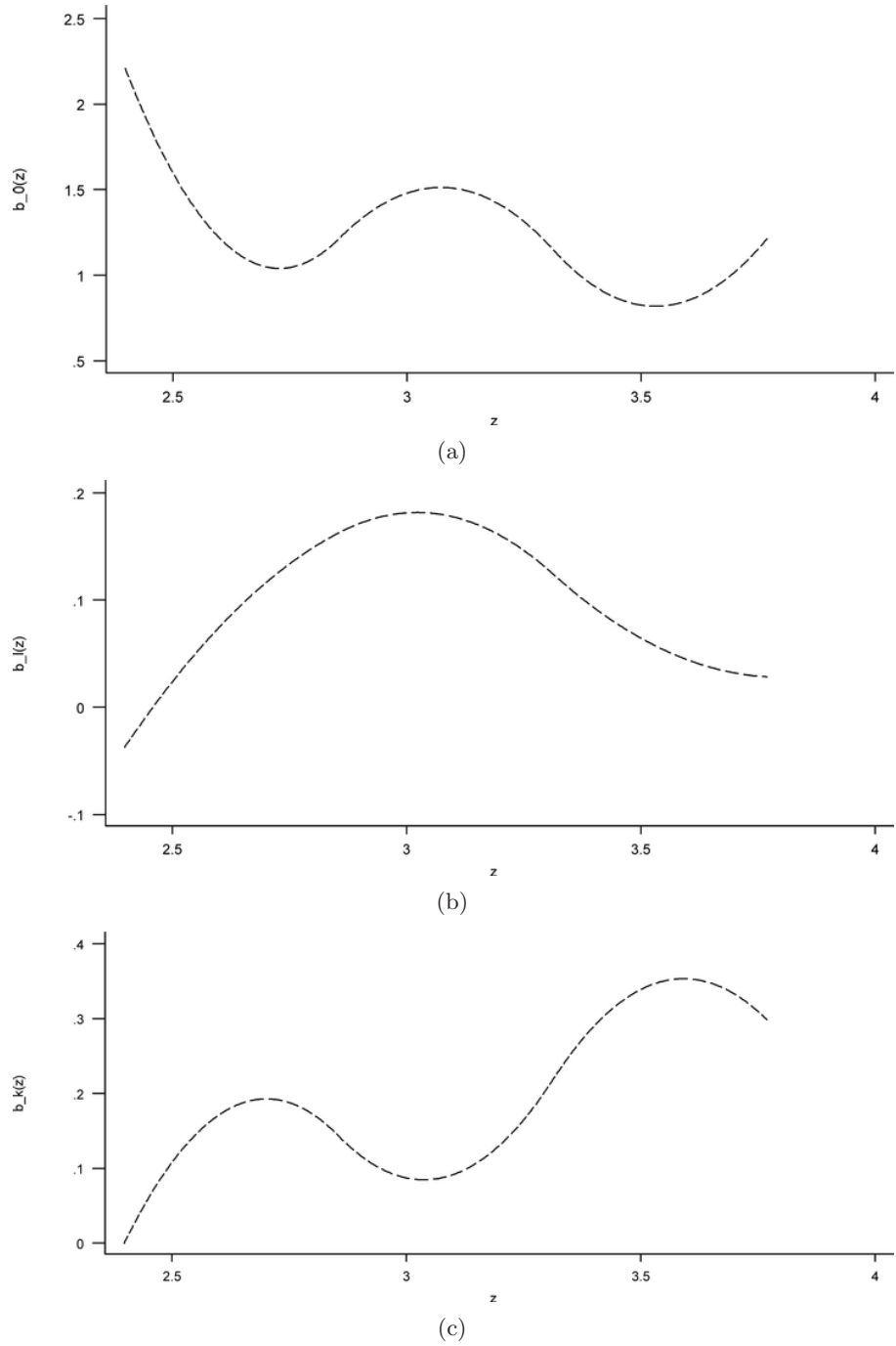

FIG. 3. (a) $b_0(z)$ *(spline)*. (b) $b_1(z)$ *(spline)*. (c) $b_k(z)$ *(spline)*.

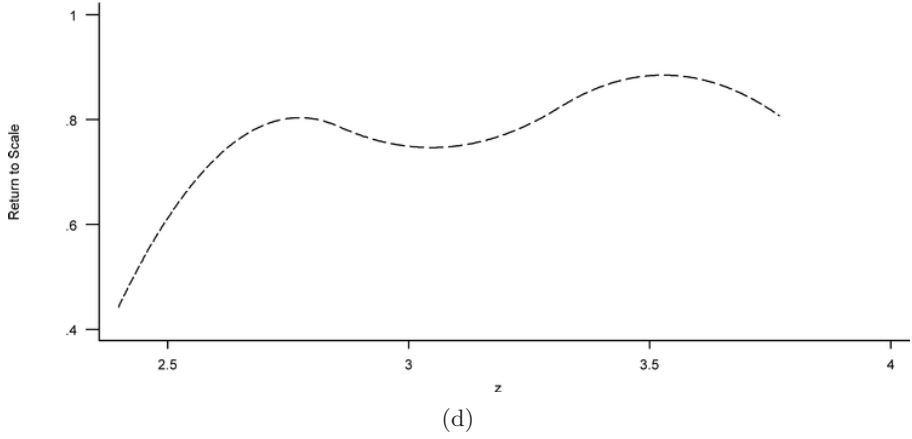

Fig. 3 (continued). (d) *Return to scale (spline)*.

semiparametric efficient. An efficient estimator can be obtained by dividing each term in (5) by $\sigma_i = \sqrt{\sigma^2(v_i)}$:

$$(24) \qquad \frac{Y_i}{\sigma_i} = \frac{w_i'}{\sigma_i}\gamma + \frac{x_i'\beta(z_i)}{\sigma_i} + \frac{u_i}{\sigma_i}.$$

We estimate $(\gamma', \beta(z_i)')$ by the least squares regression of $Y_i/\sigma_i$ on $(w_i/\sigma_i, p^K(x_i, z_i)'/\sigma_i)$. The transformed error $u_i/\sigma_i$ becomes conditional homoskedastic. Under the assumption that $0 < \eta_1 \leq \inf_v \sigma^2(v) \leq \sup_v \sigma^2(v) \leq \eta_2 < \infty$ for some positive constants $\eta_1 < \eta_2$, by the same arguments as in the proof of Theorem 2.1, one can show that

$$\sqrt{n}(\tilde{\gamma} - \gamma) \to N(0, J_0^{-1} A_0 J_0^{-1}) = N(0, J_0^{-1}) \qquad \text{in distribution,}$$

where

$$(25) \qquad J_0 = \inf_{\xi \in \mathcal{G}} E\{[w_i - x_i'\xi(z_i)][w_i - x_i'\xi(z_i)]'/\sigma^2(v_i)\}$$

and

$$A_0 = \inf_{\xi \in \mathcal{G}} E\{[w_i - x_i'\xi(z_i)][w_i - x_i'\xi(z_i)]' u_i^2/\sigma^4(v_i)\}$$
$$= \inf_{\xi \in \mathcal{G}} E\{[w_i - x_i'\xi(z_i)][w_i - x_i'\xi(z_i)]'/\sigma^2(v_i)\} = J_0.$$

Therefore, by the result of Chamberlain (1992), we know that $\tilde{\gamma}$ is semiparametrically efficient. Note that if we let $a(x, z) = x'\xi(z) \in \mathcal{G}$ denote the solution of the minimization problem of (25), that is, $E\{[w_i - a(x_i, z_i)][w_i - a(x_i, z_i)]'/\sigma^2(v_i)\} = \inf_{\xi \in \mathcal{G}} E\{[w_i - x_i'\xi(z_i)][w_i - x_i'\xi(z_i)]'/\sigma^2(v_i)\}$, then $a(x, z)$, in general, differs from $m(x, z) = E_\mathcal{G}(w_i)$ defined in (16) because of the weighting function $1/\sigma^2(v_i)$.

It is unlikely that $\sigma^2(v_i)$ is known in practice. Let $\hat{\sigma}^2(v_i)$ denote a generic nonparametric estimator of $\sigma^2(v_i)$, and write $\hat{\sigma}_i = \sqrt{\hat{\sigma}^2(v_i)}$. Then one can obtain feasible estimators for $\gamma$ and $\beta(z)$ by regressing $Y_i/\hat{\sigma}_i$ on $[w_i'/\hat{\sigma}_i, p^K(x_i, z_i)/\hat{\sigma}_i]$. The resulting estimator of $\gamma$ will be semiparametric efficient provided that $\hat{\sigma}(v)$ converges to $\sigma(v)$ uniformly with a certain rate for all $v$ in the compact support of $v$.

For the kernel-based profile likelihood approach, it is more difficult to obtain efficient estimation when the error is conditional heteroskedastic. Recall that $E_{\mathcal{G}}(A_i)$ denotes the projection of $A_i$ on the varying coefficient functional space $\mathcal{G}$. From (5) we have

$$(26) \qquad y_i - E_{\mathcal{G}}(y_i) = (w_i - E_{\mathcal{G}}(w_i))'\gamma + u_i.$$

Dividing each term in (26) by $\sigma_i$, we get

$$(27) \qquad \frac{y_i - E_{\mathcal{G}}(y_i)}{\sigma_i} = \frac{(w_i - E_{\mathcal{G}}(w_i))'}{\sigma_i}\gamma + \frac{u_i}{\sigma_i}.$$

Let $\bar{\gamma}$ denote the least squares estimator of $\gamma$ based on (27). By the Lindeberg central limit theorem, we have

$$(28) \qquad \sqrt{n}(\bar{\gamma} - \gamma) \to N(0, \{E[(w_i - E_{\mathcal{G}}(w_i))(w_i - E_{\mathcal{G}}(w_i))'/\sigma_i^2]\}^{-1})$$

in distribution.

However, $\bar{\gamma}$ is not semiparametrically efficient because

$$E[(w_i - E_{\mathcal{G}}(w_i))(w_i - E_{\mathcal{G}}(w_i))'/\sigma_i^2]$$
$$\neq \inf_{g \in \mathcal{G}} E[(w_i - g(x_i, z_i))(w_i - g(x_i, z_i))'/\sigma^2(v_i)]$$

due to the weight function $1/\sigma_i^2$. [$E_{\mathcal{G}}(w_i)$ is defined as the (un-weighted) projection of $w_i$ on the varying coefficient functional space $\mathcal{G}$. It differs from the weighted projection in general.] We conjecture that some iterative procedure (similar to the backfitting algorithm) is needed in order to obtain an efficient kernel-based estimator for $\gamma$ when the error is conditional heteroskedastic.

## APPENDIX

Throughout this Appendix, $C$ denotes a generic positive constant that may be different in different uses, $\sum_i = \sum_{i=1}^n$. The norm $\|\cdot\|$ for a matrix $A$ is defined by $\|A\| = [\text{tr}(A'A)]^{1/2}$. Also, when $A$ is a matrix and $a_n$ is a positive sequence depending on $n$, $A = O_p(a_n)$ [or $o_p(a_n)$] means that each element of $A$ is $O_p(a_n)$ [or $o_p(a_n)$]. Also, when we write $A \leq C$ for a constant scalar $C$, it means that each element of $A$ is less or equal to $C$.

PROOF OF THEOREM 2.1. Recall that $\theta(x_i, z_i) = E[w_i|x_i, z_i]$, $m(x_i, z_i) = E_{\mathcal{G}}(w_i) = E_{\mathcal{G}}(\theta(x_i, z_i))$ and $\varepsilon_i = w_i - m(z_i, x_i)$. Define $v_i = w_i - \theta(x_i, z_i)$ and $\eta_i = \theta(z_i, x_i) - m(x_i, z_i)$. We will use the following short-hand notation: $\theta_i = \theta(x_i, z_i)$, $g_i = x_i'\beta(z_i)$ and $m_i = m(x_i, z_i)$. Hence, $v_i = w_i - \theta_i$, $\varepsilon_i = \theta_i + v_i - m_i$, $\eta_i = \theta_i - m_i$. Finally, the variables without subscript represent matrices, for example, $\theta = (\theta_1, \ldots, \theta_n)'$ is of dimension $n \times 1$.

Also recall that for any matrix $A$ with $n$ rows, we define $\widetilde{A} = P(P'P)^- P'A$ [$P$ is defined below (6)]. Applying this definition to $\theta, m, g, \eta, u, v$, we get $\tilde\theta, \tilde m, \tilde g, \tilde\eta, \tilde u, \tilde v$.

Since $w_i = \theta_i + v_i$ and $\theta_i = m_i + \eta_i$, we get $w_i = \eta_i + v_i + m_i$ and $\tilde w_i = \tilde\eta_i + \tilde v_i + \tilde m_i$. In matrix notation,

$$W = \eta + v + m \quad \text{and} \quad \widetilde{W} = \tilde\eta + \tilde v + \tilde m.$$

Therefore, we have

(29) $$W - \widetilde{W} = \eta + v + (m - \tilde m) - \tilde v - \tilde\eta.$$

For scalars or column vectors $A_i$ and $B_i$, we define $S_{A,B} = n^{-1}\sum_i A_i B_i'$ and $S_A = S_{A,A}$. We also define the scalar function $\overline{S}_A = n^{-1}\sum_i A_i' A_i$, which is the sum of the diagonal elements of $S_A$. Using $ab \leq (a^2 + b^2)/2$, it is easy to see that each element of $S_{A,B}$ is less or equal to $\overline{S}_A + \overline{S}_B$. When we evaluate the probability order of $S_{A,B}$, we often write $S_{A,B} \leq \overline{S}_A + \overline{S}_B$. The scalar bound $\overline{S}_A + \overline{S}_B$ bounds each of the elements in $S_{A,B}$. Therefore, if $\overline{S}_A + \overline{S}_B = O_p(a_n)$ (for some sequence $a_n$), then each element of $S_{A,B}$ is at most $O_p(a_n)$, which implies that $S_{A,B} = O_p(a_n)$. Similarly, using the Cauchy–Schwarz inequality, we have $S_{A,B} \leq (\overline{S}_A \overline{S}_B)^{1/2}$. Here again, the scalar bounds all the elements in $S_{A,B}$.

Note that if $S^{-1}_{W-\widetilde{W}}$ exists, then, from (10) and (11), we get

(30) $$\sqrt{n}(\hat\gamma - \gamma) = \left[n^{-1}\sum_i (w_i - \tilde w_i)(w_i - \tilde w_i)'\right]^{-1}$$
$$\times \sqrt{n}\left\{n^{-1}\sum_i (w_i - \tilde w_i)(g_i - \tilde g_i + u_i - \tilde u_i)\right\}$$
$$= S^{-1}_{W-\widetilde{W}}\sqrt{n} S_{W-\widetilde{W}, g-\tilde g + u - \tilde u},$$

where $g_i = x_i'\beta(z_i)$.

For the first part of the theorem, we will prove the following: (i) $S_{W-\widetilde{W}} = \Phi + o_p(1)$, (ii) $S_{W-\widetilde{W}, g-\tilde g} = o_p(n^{-1/2})$, (iii) $S_{W-\widetilde{W}, \tilde u} = o_p(n^{-1/2})$ and (iv) $\sqrt{n} S_{W-\widetilde{W}, u} \to N(0, \Omega)$ in distribution.

PROOF OF (i). For a matrix $A$ and scalar sequence $a_n$, $A = O_p(a_n)$ ($o_p(a_n)$) means that each element of $A$ has an order of $O_p(a_n)$ ($o_p(a_n)$). Using (29), we have

(31) $S_{W-\widetilde{W}} = S_{\eta+v+(m-\tilde m)-\tilde v-\tilde\eta} = S_{\eta+v} + S_{(m-\tilde m)-\tilde v-\tilde\eta} + 2S_{\eta+v,(m-\tilde m)-\tilde v-\tilde\eta}.$

The first term $S_{\eta+v} = \frac{1}{n}\sum_i (\eta_i + v_i)(\eta_i + v_i)' = \frac{1}{n}\sum_i \varepsilon_i \varepsilon_i' = \Phi + o_p(1)$ by virtue of the law of large numbers.

The second term $S_{(m-\tilde{m})-\tilde{v}-\tilde{\eta}} \leq 3(\overline{S}_{(m-\tilde{m})} + \overline{S}_{\tilde{v}} + \overline{S}_{\tilde{\eta}}) = o_p(1)$ by Lemmas A.3, A.4(i) and A.5, stated and proved at the end of this Appendix.

The last term $S_{\eta+v,(m-\tilde{m})-\tilde{v}-\tilde{\eta}} \leq \{\overline{S}_{\eta+v}\overline{S}_{(m-\tilde{m})-\tilde{v}-\tilde{\eta}}\}^{1/2} = (O_p(1)o_p(1))^{1/2} = o_p(1)$ by the preceding results, where for an $m \times m$ matrix $A$, $\mathrm{Diag}(A)$ is an $m \times 1$ matrix with the diagonal elements of $A$, and $A^{1/2}$ has the same dimension as $A$ by taking the square root for each element of $A$. $\square$

PROOF OF (ii). Using (29), we have

$$S_{W-\widetilde{W},g-\tilde{g}} = S_{\eta+v+(m-\tilde{m})-\tilde{v}-\tilde{\eta},g-\tilde{g}} \tag{32}$$
$$= S_{\eta+v,g-\tilde{g}} + S_{m-\tilde{m},g-\tilde{g}} - S_{\tilde{v},g-\tilde{g}} - S_{\tilde{\eta},g-\tilde{g}}.$$

For the first term, by noting that $\eta_i + v_i$ is orthogonal to the varying coefficient functional space $\mathcal{G}$, and $g_i - \tilde{g}_i$ belong to $\mathcal{G}$, we have using Lemma A.3, $E[\|S_{\eta+v,g-\tilde{g}}\|^2] = n^{-2}\sum_{i=1}^n E[(\eta_i + v_i)(\eta_i + v_i)'(g_i - \tilde{g}_i)^2] \leq Cn^{-1}(\sum_{l=1}^d k_l^{2\delta_l}) \times E[\|\eta_1 + v_1\|^2] = O(n^{-1}\sum_{l=1}^d k_l^{2\delta_l}) = o(n^{-1})$, which implies that $S_{\eta+v,g-\tilde{g}} = O_p(n^{-1/2}\sum_{l=1}^d k_l^{-\delta_l})$.

The second term $S_{m-\tilde{m},g-\tilde{g}} \leq (\overline{S}_{m-\tilde{m}}\overline{S}_{g-\tilde{g}})^{1/2} = O_p(\sum_{l=1}^d k_l^{-2\delta_l})$ by Lemma A.3.

The third term $S_{\tilde{v},g-\tilde{g}} \leq (\overline{S}_{\tilde{v}}\overline{S}_{g-\tilde{g}})^{1/2} = O_p((K/n)^{1/2})O_p(\sum_{l=1}^d k_l^{-\delta_l})$ by Lemmas A.3 and A.4(i). The last term $S_{\tilde{\eta},g-\tilde{g}} \leq (\overline{S}_{\tilde{\eta}}\overline{S}_{g-\tilde{g}})^{1/2} = O_p((k/n)^{1/2}) \times O_p(\sum_{l=1}^d k_l^{-\delta_l})$ by Lemmas A.3 and A.5.

Combining the above four terms we have $S_{W-\widetilde{W},g-\tilde{g}} = O_p((n^{-1/2} + (K/n)^{1/2})(\sum_{l=1}^d k_l^{-\delta_l}) + \sum_{l=1}^d k_l^{-2\delta_l}) = o_p(n^{-1/2})$ by Assumption 2.3. $\square$

PROOF OF (iii). Using (29), we have

$$S_{W-\widetilde{W},\tilde{u}} = S_{\eta+v+(m-\tilde{m})-\tilde{v}-\tilde{\eta},\tilde{u}} = S_{\eta+v,\tilde{u}} + S_{m-\tilde{m},\tilde{u}} - S_{\tilde{v},\tilde{u}} - S_{\tilde{\eta},\tilde{u}}. \tag{33}$$

The first term $S_{\eta+v,\tilde{u}} \leq (\overline{S}_{\eta+v}\overline{S}_{\tilde{u}})^{1/2} = O_p(K/n)$ by Lemma A.4(ii). The second term $S_{m-\tilde{m},\tilde{u}} \leq (\overline{S}_{m-\tilde{m}}\overline{S}_{\tilde{u}})^{1/2} = O_p(\sum_{l=1}^d k_l^{-\delta_l})O_p(\sqrt{K}/\sqrt{n})$ by Lemmas A.3 and A.4(ii).

The third term $S_{\tilde{v},\tilde{u}} \leq (\overline{S}_{\tilde{v}}\overline{S}_{\tilde{u}})^{1/2} = O_p(K/n)$ by Lemma A.4(i), (ii). The last term $S_{\tilde{\eta},\tilde{u}} \leq (\overline{S}_{\tilde{\eta}}\overline{S}_{\tilde{u}})^{1/2} = O_p(K/n)$ by Lemmas A.4(ii) and A.5.

Combining all four terms, we get $S_{W-\widetilde{W},\tilde{u}} = O_p(K/n + n^{-1/2}\sum_{l=1}^d k_l^{-\delta_l}) = o_p(n^{-1/2})$ by Assumption 2.3. $\square$

PROOF OF (iv). Using (29), we have

$$\sqrt{n}S_{W-\widetilde{W},u} = \sqrt{n}S_{\eta+v+(m-\tilde{m})-\tilde{v}-\tilde{\eta},u} \tag{34}$$
$$= \sqrt{n}S_{\eta+v,u} + \sqrt{n}(S_{m-\tilde{m},u} - S_{\tilde{v},u} - S_{\tilde{\eta},u}).$$

The first term $\sqrt{n}S_{\eta+v,u} = \sqrt{n}\sum_{i=1}^n(\eta_i+v_i)u_i = \sqrt{n}\sum_{i=1}^n\varepsilon_i u_i \to N(0,\Omega)$ in distribution by the Lindeberg–Feller central limit theorem.

The second term $E[S^2_{m-\tilde{m},u}|X,Z] = \frac{1}{n^2}\operatorname{tr}\{(m-\tilde{m})(m-\tilde{m})'E[uu'|X,Z]\} \leq (C/n)\operatorname{tr}[(m-\tilde{m})'(m-\tilde{m})/n] = (C/n)S_{m-\tilde{m}} = o_p(n^{-1})$ by Lemma A.3. Hence, $S_{m-\tilde{m},u} = o_p(n^{-1/2})$.

The third term $E[S^2_{\tilde{v},u}|X,Z] = \frac{1}{n^2}\operatorname{tr}(P(P'P)^{-1}P'vv'P(P'P)^{-1}P'E[uu'|X,Z]) \leq (C/n^2)\operatorname{tr}[P(P'P)^{-1}P'vv'P(P'P)^{-1}P'] = (C/n)\operatorname{tr}(\tilde{v}\tilde{v}'/n) = (C/n)S_{\tilde{v}} = o_p(n^{-1})$ by Lemma A.4(i). Hence, $S_{\tilde{v},u} = o_p(n^{-1/2})$.

The last term $S_{\tilde{\eta},u} = o_p(n^{-1/2})$ by the same proof as $S_{\tilde{v},u} = o_p(n^{-1/2})$ by citing Lemma A.5, rather than citing Lemma A.4(i). □

Combining proofs of (i)–(iv) with (30), we conclude that $\sqrt{n}(\hat{\gamma}-\gamma) \to N(0, \Phi^{-1}\Omega\Phi^{-1})$ in distribution.

For the second part of the theorem, we need to show that $\hat{\Sigma} = \Sigma + o_p(1)$, where $\hat{\Sigma} = \hat{\Phi}^{-1}\hat{\Omega}\hat{\Phi}^{-1}$. But $\hat{\Phi} = S_{W-\widetilde{W}} = \Phi + o_p(1)$ is proved in the proof of (i) above. By a similar argument, it is easy to show that $\hat{\Omega} = \Omega + o_p(1)$. Therefore, $\hat{\Sigma} = \Sigma + o_p(1)$. □

PROOF OF THEOREM 2.2. We will prove Theorem 2.2 by replacing $\hat{\beta}(z)$ and $\beta(z)$ by $\hat{g}(x,z) = x'\hat{\beta}(z)$ and $g(x,z) = x'\beta(z)$, respectively, because $|\hat{g}(x,z) - g(x,z)|^2 = |x'(\hat{\beta}(z) - \beta(z))|^2 \leq d\sum_{l=1}^d x_l^2(\hat{\beta}_l(z) - \beta_l(z))^2$, which has the same order as $\|\hat{\beta}(z) - \beta(z)\|^2$ under the bounded support assumption. Hence, the rate of convergence for $\hat{g}(x,z) - g(x,z)$ is the same as that of $\hat{\beta}(z) - \beta(z)$.

The proof is similar to the proof of Theorem 1 in Newey (1997). Define an indicator function $\mathbf{1}_n$ which equals 1 if $(P'P)$ is nonsingular and 0 otherwise. We first find the convergence rate of $\mathbf{1}_n\|\hat{\alpha} - \alpha\|$. By (12) and (7), and if $(P'P)^{-1}$ exists, we have

$$\hat{\alpha} = (P'P)^{-1}P'(Y - W\hat{\gamma})$$
$$= (P'P)^{-1}P'(Y - W\gamma - W(\hat{\gamma}-\gamma))$$
(35) $$= (P'P)^{-1}P'(P\alpha + (G - P\alpha) + u - W(\hat{\gamma}-\gamma))$$
$$= \alpha + (P'P/n)^{-1}P'(G-P\alpha)/n + (P'P/n)^{-1}P'u/n$$
$$\quad - (P'P/n)^{-1}P'W(\hat{\gamma}-\gamma)/n.$$

Hence,

(36)
$$\mathbf{1}_n\|\hat{\alpha}-\alpha\| \leq \mathbf{1}_n\|(P'P/n)^{-1}P'(G-P\alpha)/n\|$$
$$+ \mathbf{1}_n\|(P'P/n)^{-1}P'u/n\|$$
$$+ \mathbf{1}_n\|(P'P/n)^{-1}P'W(\hat{\gamma}-\gamma)/n\|.$$

The first term $\mathbf{1}_n\|(P'P/n)^{-1}P'(G-P\alpha)/n\| = O_p(\sum_{l=1}^{d} k_l^{-\delta_l})$ by Lemma A.2.

The second term

$$E[\mathbf{1}_n\|(P'P/n)^{-1}P'u/n\| | X, Z]$$
$$= \mathbf{1}_n E[((u'P/n)(P'P/n)^{-1}(P'P/n)^{-1}(P'u/n))^{1/2}|X,Z]$$
$$\leq O_p(1)\mathbf{1}_n \operatorname{tr}(P(P'P)^{-1}P'E[uu'|X,Z]/n)^{1/2}$$
$$\leq O_p(1)\mathbf{1}_n C\sqrt{K}/\sqrt{n}$$

by Lemma A.1 and Assumption 2.1. Hence, $\mathbf{1}_n\|(P'P/n)^{-1}P'u/n\| = O_p(\sqrt{K}/\sqrt{n})$.

As for the last term, note that $W = \eta + v + m = \varepsilon + m$ and $\hat{\gamma} - \gamma = O_p(n^{-1/2})$ by Theorem 2.1. Therefore,

$$E[\mathbf{1}_n\|(P'P/n)^{-1}P'W/n\| | X, Z]$$
$$= \mathbf{1}_n E[\|(P'P/n)^{-1}P'(\varepsilon+m)/n\| | X, Z]$$
$$\leq \mathbf{1}_n E[\|(P'P/n)^{-1}P'\varepsilon/n\| | X, Z] + \mathbf{1}_n E[\|(P'P/n)^{-1}P'm/n\| | X, Z].$$

Also,

$$\mathbf{1}_n E[\|(P'P/n)^{-1}P'\varepsilon/n\| | X, Z]$$
$$= \mathbf{1}_n E[\|(\varepsilon'P/n)(P'P/n)^{-1}(P'P/n)^{-1}(P'\varepsilon/n)\| | X, Z]$$
$$\leq O_p(1)\mathbf{1}_n \operatorname{tr}(P(P'P)^{-1}P'E[\varepsilon\varepsilon'|X,Z]/n)^{1/2}$$
$$\leq O_p(1)\mathbf{1}_n C\sqrt{K}/\sqrt{n}$$

by Lemma A.1 as in the proof of Theorem 2.1. Hence, $\mathbf{1}_n|(P'P/n)^{-1}P'\varepsilon/n| = O_p(\sqrt{K}/\sqrt{n}) = o_p(1)$.

$\mathbf{1}_n\|(P'P)^{-1}P'm\| = \mathbf{1}_n\|(P'P/n)^{-1}P'm/n\| = O_p(1)$ by Lemma A.2.

Combining the above results, also noting that $\mathbf{1}_n \to 1$ almost surely, we have

(37) $$\|\hat{\alpha} - \alpha\| = O_p\left(\sum_{l=1}^{d} k_l^{-\delta_l} + \sqrt{K}/\sqrt{n}\right).$$

To prove part (i) of Theorem 2.2, using (37) and Assumption 2.3, and also noting that $\hat{g}(x,z) = x'\hat{\beta}(z) = p^K(x,z)'\hat{\alpha}$, we have

$$\sup_{(x,z)\in\mathcal{S}} |\hat{g}(x,z) - g(x,z)| \leq \sup_{(x,z)\in\mathcal{S}} |p^K(x,z)'(\hat{\alpha}-\alpha)| + |p^K(x,z)'\alpha - g(x,z)|$$

$$\leq \zeta_0(K)\|\hat{\alpha}-\alpha\| + O\left(\sum_{l=1}^{d} k_l^{-\delta_l}\right)$$

$$= O_p\bigg(\zeta_0(K)\bigg(\sum_{l=1}^{d} k_l^{-\delta_l} + \sqrt{K}/\sqrt{n}\bigg)\bigg).$$

Proofs for (ii) and (iii) are similar, and we only prove (ii),

$$n^{-1}\sum_{i=1}^{n}[\hat{g}(x_i, z_i) - g(z_i, z_i)]^2$$
$$= n^{-1}\|P\hat{\alpha} - G\|^2$$
$$\leq 2n^{-1}\{\|P(\hat{\alpha} - \alpha)\|^2 + \|P\alpha - G\|^2\}$$
$$= 2(\hat{\alpha} - \alpha)'(P'P/n)(\hat{\alpha} - \alpha) + 2\sup_{(x,z)\in\mathcal{S}}[p^K(x,z)\alpha - g(x,z)]^2$$
$$= O_p\bigg(K/n + \sum_{l=1}^{d} k_l^{-2\delta_l}\bigg)$$

by (37), Lemma A.1 and Assumption 2.3(i). Thus, we have proved Theorem 2.2.

□

We now present some lemmas that are used in the proofs of Theorems 2.1 and 2.2. We will omit the indicator function $\mathbf{1}_n$ below since $\text{Prob}(\mathbf{1}_n = 1) \to 1$ almost surely. Following the arguments in Newey (1997), we can assume without loss of generality that $B = I$ ($B$ is defined in Assumption 2.2). Hence, $P^K(X, Z) = p^K(X, Z)$, and $Q = E[p^K(x_i, z_i)p^K(x_i, z_i)'] = I$ ($I$ is an identity matrix of dimension $K$); see Newey (1997) for the reasons and more discussion of these issues. Recall that $p^K(x, z)$ is a $K \times 1$ matrix and rewrite each component of this matrix as $p^K(x, z) = (p_{1K}(x, z), \ldots, p_{KK}(x, z))'$.

LEMMA A.1. $\|\hat{Q} - I\| = O_p(\zeta_0(K)\sqrt{K}/\sqrt{n}) = o_p(1)$, where $\hat{Q} = P'P/n$.

PROOF. This is Theorem 1 in Newey (1997). □

LEMMA A.2. $\|\tilde{\alpha}_f - \alpha_f\| = O_p(\sum_{l=1}^{d} k_l^{-\delta_l})$, where $\tilde{\alpha}_f = (P'P)^{-1}P'f$, $\alpha_f$ satisfies Assumption 2.3 and $f = G$ or $f = m$.

PROOF. By Lemma A.1, Assumption 2.3 and the fact that $P(P'P)^{-1}P'$ is idempotent,

$$\|\tilde{\alpha}_f - \alpha_f\| = \|(P'P)^{-1}P'(f - P\alpha_f)\|$$
$$= \|(f - P\alpha_f)'P(P'P)^{-1}\hat{Q}P'(f - P\alpha_f)/n\|^{1/2}$$
$$\leq O_p(1)\|(f - P\alpha_f)'P(P'P)^{-1}P'(f - P\alpha_f)/n\|^{1/2}$$

$$\leq O_p(1)\|(f - P\alpha_f)'(f - P\alpha_f)/n\|^{1/2} = O_p\bigg(\sum_{l=1}^d k_l^{-\delta_l}\bigg). \qquad \square$$

LEMMA A.3. $S_{f-\tilde{f}} = O_p(\sum_{l=1}^d k_l^{-2\delta_l})$, where $f = G$ or $f = m$.

PROOF. Note that $\tilde{f} = P\tilde{\alpha}_f$. By Assumption 2.3 and Lemmas A.1 and A.2,

$$S_{f-\tilde{f}} = \frac{1}{n}|f - \tilde{f}|^2 \leq \frac{1}{n}(|f - P\alpha_f|^2 + |P(\alpha_f - \tilde{\alpha}_f)|^2)$$

$$= O\bigg(\sum_{l=1}^d k_l^{-2\delta_l}\bigg) + (\alpha_f - \tilde{\alpha}_f)'(P'P/n)(\alpha_f - \tilde{\alpha}_f)$$

$$\leq O\bigg(\sum_{l=1}^d k_l^{-2\delta_l}\bigg) + O_p(1)|\alpha_f - \tilde{\alpha}_f|^2 = O_p\bigg(\sum_{l=1}^d k_l^{-2\delta_l}\bigg). \qquad \square$$

LEMMA A.4. (i) $S_{\tilde{v}} = O_p(K/n)$, (ii) $S_{\tilde{u}} = O_p(K/n)$.

PROOF. (i) This proof is similar to the proof of Theorem 1 of Newey (1997),

$$E[S_{\tilde{v}}|X,Z] = \frac{1}{n}E[v'P(P'P)^{-1}P'v|X,Z]$$

$$= \frac{1}{n}E[\text{tr}(P(P'P)^{-1}P'E[vv'|X,Z])]$$

$$\leq \frac{C}{n}\text{tr}(P(P'P)^{-1}P') = C\bigg(\frac{K}{n}\bigg).$$

Hence, $S_{\tilde{v}} = O_p(K)$.
(ii) follows as in the proof of Lemma A.4(i). $\square$

LEMMA A.5. $S_{\tilde{\eta}} = O_p(K/n)$.

PROOF. First we show that $(P'\eta/n) = O_p(\sqrt{K}/\sqrt{n})$. Recall that $\theta(x_i, z_i) = E(w_i|x_i, z_i)$ and $\eta_i = \theta(x_i, z_i) - E_{\mathcal{G}}[\theta(x_i, z_i)]$. Note that $p^K(x_i, z_i) \in \mathcal{G}$ and $E_{\mathcal{G}}(\eta_i) = 0$ (i.e., $\eta \perp \mathcal{G}$). Hence, $E\|P'\eta/n\|^2 = n^{-2}\sum_i E[p^K(x_i)'\|\eta_i\|^2 p^K(x_i)] \leq \frac{C}{n}E[p^K(X_i)'p^K(x_i)] = \frac{C}{n}\text{tr}\{E[p^K(X_i)p^K(x_i)']\} = (CK/n) = O(K/n)$, which implies that $(P'\eta/n) = O_p(\sqrt{K}/\sqrt{n})$.

Thus, $S_{\tilde{\eta}} = n^{-1}\tilde{\eta}'\eta = (\eta'P/n)(P'P/n)^{-1}(P'\eta/n) = O_p(K/n)O_p(1) = O_p(K/n)$ by Lemma A.1 and the fact that $P'\eta/n = O_p(\sqrt{K}/\sqrt{n})$ as shown above. $\square$

**Acknowledgments.** We are grateful for the insightful comments from two referees and an Associate Editor which greatly improved the paper. We would also like to thank Dong Li for his help in the empirical work of this paper.

## REFERENCES

AI, C. and CHEN, X. (2003). Efficient estimation of models with conditional moment restrictions containing unknown functions. *Econometrica* **71** 1795–1843. MR2015420

ANDREWS, D. W. K. (1991). Asymptotic normality of series estimators for nonparametric and semiparametric regression models. *Econometrica* **59** 307–345. MR1097531

BICKEL, P. J., KLAASSEN, C. A. J., RITOV, Y. and WELLNER, J. A. (1993). *Efficient and Adaptive Inference for Semiparametric Models*. Johns Hopkins Univ. Press. MR1245941

BICKEL, P. J. and KWON, J. (2002). Inference for semiparametric models: Some current frontiers and an answer (with discussion). *Statist. Sinica* **11** 863–960. MR1867326

CAI, Z., FAN, J. and LI, R. (2000). Efficient estimation and inferences for varying-coefficient models. *J. Amer. Statist. Assoc.* **95** 888–902. MR1804446

CAI, Z., FAN, J. and YAO, Q. (2000). Functional-coefficient regression models for nonlinear time series. *J. Amer. Statist. Assoc.* **95** 941–956. MR1804449

CARROLL, R. J., FAN, J., GIJBELS, I. and WAND, M. P. (1997). Generalized partially linear single-index models. *J. Amer. Statist. Assoc.* **92** 477–489. MR1467842

CHAMBERLAIN, G. (1992). Efficiency bounds for semiparametric regression. *Econometrica* **60** 567–596. MR1162999

CHEN, R. and TSAY, R. S. (1993). Functional-coefficient autoregressive models. *J. Amer. Statist. Assoc.* **88** 298–308. MR1212492

CRAVEN, P. and WAHBA, G. (1979). Smoothing noisy data with spline functions: Estimating the correct degree of smoothing by generalized cross-validation. *Numer. Math.* **31** 377–403. MR516581

FAN, J. and HUANG, L.-S. (2001). Goodness-of-fit tests for parametric regression models. *J. Amer. Statist. Assoc.* **96** 640–652. MR1946431

FAN, J. and HUANG, T. (2002). Profile likelihood inferences on semiparametric varying-coefficient partially linear models. Unpublished manuscript.

FAN, J., YAO, Q. and CAI, Z. (2003). Adaptive varying-coefficient linear models. *J. R. Stat. Soc. Ser. B Stat. Methodol.* **65** 57–80. MR1959093

FAN, J. and ZHANG, W. (1999). Statistical estimation in varying coefficient models. *Ann. Statist.* **27** 1491–1518. MR1742497

GREEN, P. J. and SILVERMAN, B. W. (1994). *Nonparametric Regression and Generalized Linear Models*: *A Roughness Penalty Approach*. Chapman and Hall, London. MR1270012

HÄRDLE, W., LIANG, H. and GAO, J. (2000). *Partially Linear Models*. Physica-Verlag, Heidelberg. MR1787637

HART, J. (1997). *Nonparametric Smoothing and Lack-of-fit Tests*. Springer, New York. MR1461272

HASTIE, T. and TIBSHIRANI, R. (1993). Varying coefficient models (with discussion). *J. Roy. Statist. Soc. Ser. B* **55** 757–796. MR1229881

HOOVER, D. R., RICE, J. A., WU, C. O. and YANG, L.-P. (1998). Nonparametric smoothing estimates of time-varying coefficient models with longitudinal data. *Biometrika* **85** 809–822. MR1666699

HUANG, J. Z. (1998). Projection estimation in multiple regression with application to functional ANOVA models. *Ann. Statist.* **26** 242–272. MR1611780

Huang, J. Z. (2003). Local asymptotics for polynomial spline regression. *Ann. Statist.* **31** 1600–1635. MR2012827

Huang, J. Z., Wu, C. O. and Zhou, L. (2002). Varying coefficient models and basis function approximations for the analysis of repeated measurements. *Biometrika* **89** 111–128. MR1888349

Huang, J. Z., Wu, C. O. and Zhou, L. (2004). Polynomial spline estimation and inference for varying coefficient models with longitudinal data. *Statist. Sinica* **14** 763–788. MR2087972

Li, K. C. (1987). Asymptotic optimality for $C_p$, $C_L$, cross-validation and generalized cross-validation: Discrete index set. *Ann. Statist.* **15** 958–975. MR902239

Li, Q., Huang, C. J., Li, D. and Fu, T.-T. (2002). Semiparametric smooth coefficient models. *J. Bus. Econom. Statist.* **20** 412–422. MR1939909

Lorentz, G. G. (1966). *Approximation of Functions*. Holt, Rinehart and Winston, New York. MR213785

Mallows, C. L. (1973). Some comments on $C_p$. *Technometrics* **15** 661–675.

Newey, W. K. (1997). Convergence rates and asymptotic normality for series estimators. *J. Econometrics* **79** 147–168. MR1457700

Robinson, P. M. (1988). Root-N-consistent semiparametric regression. *Econometrica* **56** 931–954. MR951762

Shen, X. (1997). On methods of sieves and penalization. *Ann. Statist.* **25** 2555–2591. MR1604416

Speckman, P. (1988). Kernel smoothing in partially linear models. *J. Roy. Statist. Soc. Ser. B* **50** 413–436. MR970977

Stock, C. J. (1989). Nonparametric policy analysis. *J. Amer. Statist. Assoc.* **89** 567–575. MR1010347

Xia, Y. C. and Li, W. K. (1999). On the estimation and testing of functional-coefficient linear models. *Statist. Sinica* **9** 735–758. MR1711643

Xia, Y. C. and Li, W. K. (2002). Asymptotic behavior of bandwidth selected by the cross-validation method for local polynomial fitting. *J. Multivariate Anal.* **83** 265–287. MR1945954

Zhang, W., Lee, S.-Y. and Song, X. (2002). Local polynomial fitting in semivarying coefficient models. *J. Multivariate Anal.* **82** 166–188. MR1918619

Zhou, S., Shen, X. and Wolfe, D. A. (1998). Local asymptotics for regression splines and confidence regions. *Ann. Statist.* **26** 1760–1782. MR1673277

I. Ahmad
Department of Statistics
University of Central Florida
Orlando, Florida 32816-2370
USA
e-mail: iahmad@mail.ucf.edu

S. Leelahanon
Faculty of Economics
Thammasat University
2 Pra Chan Road
Bangkok 10200
Thailand

Q. Li
Department of Economics
Texas A&M University
College Station, Texas 77843-4228
USA